\documentclass[reqno, 12pt]{amsart}
\def\C{{\mathbb C}}

\def\Z{{\mathbb Z}}
\def\Zp{{\mathbb Z}_p}
\def\R{{\mathbb R}}

\def\Fp{{{\mathbb F}_p}}
\def\F{{\mathbb F}}

\def\Fq{{{\mathbb F}_q}}

\newtheorem{cor}{Corollary}

\newtheorem{theorem}{Theorem}

\theoremstyle{definition}
\newtheorem{defn}{Definition}   

\newtheorem{question}{Question}
\newtheorem{convention}{Convention}

\theoremstyle{remark}
\newtheorem{rem}{Remark}

\begin{document}
\title[A local field approach to the Riemann Hypothesis]
{A local field approach to the Riemann Hypothesis}
%\subjclass{Primary 11G09}
%\keywords{Drinfeld modules, $T$-modules, characteristic $p$ $L$-series,
%$L$-zeroes, multi-valued operators}
\author{David Goss}
%\thanks{}
\address{Department of Mathematics\\The Ohio State University\\231
W.\ $18^{\rm th}$ Ave.\\Columbus, Ohio 43210}

\email{goss@math.ohio-state.edu}

\date{June 8, 2012}

\begin{abstract}
Since the seminal work of Wan, Poonen, and Sheats in the 1990's, we have been 
searching for the correct general statement of the Riemann Hypothesis (``RH'') which appears implicit
in their results. Recently, upon viewing the extension $\C/\R$ in
light of results derived for the
Carlitz module, we were led to view the RH as a statement about ramification which we
explore in this short work.
We shall see that, combined with some ideas flowing from
the proofs of Wan, Poonen, and Sheats, this ramification idea has a good deal of explanatory
power in finite characteristic. Indeed, unramified extensions of nonarchimedean local fields are cyclotomic
in nature and this fits perfectly with the best possible extension of the result of Wan 
and Sheats.
The notion that the zeroes ``lie on a line'' seems to
be the beginning of the story in finite characteristic, and we show how this fits
with having the zeroes be unramified.

\end{abstract}

\maketitle

%****************************************************************************
\section{Introduction}\label{intro} 
In 1996, in response to a query from the author, Daqing Wan \cite{wa1} 
calculated
the Newton polygon associated to the characteristic zeta function associated to $\Fp[\theta]$. What he 
discovered was completely remarkable: there is at most one (with multiplicity) zero of a given
absolute value. As an obvious and immediate corollary, one finds that the zeroes of the zeta function
(at the infinite prime) are both  simple and lie 
``on the line'' $\Fp((1/\theta))$. Clearly this is both a ``Riemann Hypothesis,'' as the zeroes lie on a ``line'', and
the expected analog of the simplicity of zeroes of the Riemann zeta function. 

Since that time, due to the labors of D.\ Thakur, B.\ Poonen, J.\ Diaz-Vargas and J.\ Sheats \cite{dv1, sh1,
ba1}, Wan's result 
has been extended
to the general polynomial ring $\Fq[\theta]$ with the same results. While we do not yet have the theory
to understand the arithmetic applications of this result,  one clearly wants to put it into an appropriate
framework for general $L$-series. We tried this originally in \cite{go2} which is now superseded by the present work
(where
we are more specific about the fields generated by the zeroes; moreover in \cite{go2} we were not careful
enough about the factorizations of zeta functions in terms of $L$-series as discussed below in Subsection
\ref{failure}).

The zeta-function of $\Fq[\theta]$ is known to have an analytic continuation as a family $\zeta(x,y)$ of entire
functions in $x^{-1}$, where $y\in \Zp$ and $\Zp$ is the $p$-adic integers, see e.g., Section 8 of \cite{go1}, \cite{bo1}
or \cite{bp1}.
In fact, one finds that {\it all}
$L$-series of Drinfeld modules defined over finite extensions $L/\Fq(\theta)$ have similar analytic 
continuations. So a first guess at a ``generalized Riemann Hypothesis,'' would be that all the zeroes of such
$L$-series lie in $\Fq((1/\theta))$. On the other hand, nonArchimedean analysis is very algebraic in nature
and every entire function is completely determined by its zeroes (unlike, obviously, the exponential function classically).
We shall see in Section \ref{failure} that this immediately leads to counterexamples of this first guess.

Recently upon viewing the extension $\C/\R$ from 
the point of view of the Carlitz module, as in Subsection \ref{view} below, we have come to some clarity on these
issues.  Because of this insight, we view the classical Riemann hypothesis
as the statement that the zeta zeroes are ``unramified''.  We will show how, in finite characteristic, this fits very nicely
into the framework initiated by Wan and this will lead to
speculation on the nature of ``most'' zeta zeroes (and perhaps all in certain
circumstances). Indeed, combined with the calculations of
Wan/Poonen/Sheats, this new insight has some serious predictive power and actually
leads to the best possible guess on the nature of
the zeroes in view of the above counterexamples.

One of the basic statements about the Riemann zeta function is, of course, its famous
functional equation. It has long been known that a functional equation of the
same $s\mapsto 1-s$ type does not appear to hold with characteristic $p$ zeta 
functions. However, in \cite{go4} we introduced the group $S_{(q)}$ of self-homeomorphisms
of $\Zp$ (essentially $S_{(q)}$ is the permutation group on $\{0,1,2...\}$ acting
on $\Zp$ by simply permuting the $q$-adic coefficients). Remarkably, there is also
{\em some} evidence that $S_{(q)}$ acts on the zeta zeroes themselves via their
expansions as elements of $\Fq((1/\theta))$. As this action in finite
characteristic is in fact linear, we will see in Subsection
\ref{con} that it also works 
well with the new optic of ``unramified'' zeroes.

We emphasize again that we are speculating here and  actually establish very little. However, it is our goal to
present a framework to view the original Riemann Hypothesis which {\it also} leads to important ideas for characteristic
$p$ valued $L$-series and which, perhaps, gives a suitable framework for future exploration.

I thank Rudy Perkins for his careful reading of earlier versions of this work.

\subsection{Notation}\label{basic}
We set $q=p^{m_0}$ where $p$ is prime and $m_0$ is a nonnegative integer, and
let $\Fq$ be the field with $q$ elements. We set $A:=\Fq[\theta]$, $k:=\Fq(\theta)$,
and $K:=\Fq((1/\theta))$. Notice that $K$ is the completion of $k$ at the normalized
absolute value associated to the infinite place of $k$. Let $\bar K$ be a fixed algebraic closure of
$K$ and let $\bar{\F}_q\subset \bar K$ be the algebraic closure of $\Fq$.

\section{$\C/\R$ from the viewpoint of the Carlitz module}\label{view}
It has always been an interesting game whether the extension $\C/\R$ should be 
viewed as (somehow) ramified
or not; of course, as $[\C\colon \R]=2$, {\em if} one views this extension as ramified, then the ramification
should be viewed as total.

One piece of evidence that points to this extension being unramified is that the value group of
the usual archimedean absolute value $\vert ?\vert$ is the same on both fields, and this is a hallmark
of an unramified extension (of nonArchimedean local fields)...

On the other hand, the theory of function fields over a finite field
contains a number of analogs of the extension
$\C/\R$. Indeed, one analog would be our chosen algebraic closure $\bar K$ of $K$; this is clearly
algebraically closed, but is not complete. So another analog of $\C$ is the field 
$\C_\infty$ given
by the completion of $\bar K$ via the canonical extension of the absolute value 
$\vert?\vert_\infty$.
The field $\C_\infty$ is then both complete {\it and} algebraically closed. However, 
it is {\em not}
a local field (i.e., it is not locally compact); so the analogy is not exact.

However, the Carlitz module (\cite{ca1}, or Section 3 of \cite{go1}) immediately presents yet another analog of $\C$ 
that we call $K_1$ and which, in our view, mandates viewing $\C/\R$ as ramified. Let
$\exp$ be the (classical) exponential function and let $\exp_C$ be the exponential function of the
Carlitz module. Recall that $2\pi i$ is, of course, the period of $\exp$ and $\C=\R(2\pi i)$. 
Early on, Carlitz gave a beautiful  formula for the ``period'' $\tilde{\xi}$ 
of $\exp_C$ which was remarkably
simplified in  \cite{at1} in the following fashion. Let $\theta_1$ be a 
fixed choice of $(q-1)$-st root of $-\theta$. Then one has
\begin{equation}\label{view1}
\tilde{\xi}=\theta_1 \theta \prod_{i=1}^\infty \left( 1-\theta^{1-q^i}\right)^{-1}\,.
\end{equation}

\begin{defn}\label{view1.1}
We set $K_1:=K(\tilde{\xi})$. \end{defn}

Obviously $K_1$ is a local field and, from the present viewpoint, is thus clearly an analog of
$\C$ (though by no means is it algebraically closed). Also clearly, $K_1/K$ is a {\it totally (and tamely)
ramified extension} with Galois group isomorphic to $\Fq^\ast=A^\ast$ (much as $\C/\R$ has Galois group
isomorphic to $\Z^\ast$).

\begin{convention}\label{view2}
From now on, we shall view $\C/\R$ as a totally ramified extension.
\end{convention}

\begin{rem}\label{view3}
Let $L$ be a nonArchimedean local field with residue characteristic $p$. It is then
universally known that the unramified extensions of $L$ are given by adjoining
$n$-th power roots of unity with $n$ not divisible by $p$.
\end{rem}

An element $\alpha$ algebraic over a local field $L$ will be called ``unramified'' if and only if
the extension $L(\alpha)/L$ is unramified.
\section{The classical Riemann Hypothesis as a ramification statement}\label{rh1}
Let $\zeta(s)=\sum_{n=1}^\infty n^{-s}$ be the Riemann zeta function. Riemann showed that $\zeta(s)$
has a meromorphic continuation to $\C$. Let $\Gamma(s)$ be Euler's gamma function and put
\begin{equation}\label{rh11}
\Xi(s):=1/2 \pi^{-s/2}s(s-1) \Gamma(s/2)\zeta(s)\,.
\end{equation}
Riemann then showed that $\Xi(s)$ is entire and satisfies the functional equation $\Xi(s)=\Xi(1-s)$.
The {\it Riemann Hypothesis} is precisely the statement that the nontrivial zeroes of $\Xi(s)$ lie
on the line $1/2+it$ where $t\in \R$. Following Riemann, we put $\tilde{\Xi}(t):=\Xi(1/2+it)$ so
that the Riemann Hypothesis is now the statement that the zeroes of $\tilde{\Xi}(t)$ are {\it real}.
\begin{convention}\label{rh12}
In keeping with Convention \ref{view2}, we rephrase the Riemann Hypothesis as the statement
that the zeroes of $\tilde{\Xi}(t)$ are {\it unramified.}
\end{convention}

\section{The RH as a ramification statement in characteristic $p$}\label{rh2}
\subsection{The results of Wan and Sheats}\label{swan}
Let $A$ etc., be as in Subsection \ref{basic}. Let $A^\prime$ be the set of monic elements of 
$A$ and formally put
\begin{equation}\label{rh21}
\zeta_A(s)=\zeta(s)=\sum_{a\in A^\prime}a^{-s}\,.
\end{equation}
It is clear, for instance, that this series will converge to an element of $K$ if $s$ is a positive integer. 
It has been shown in various
places (see, for instance, Section 8 of \cite{go1})
that $\zeta(s)$ has an analytic continuation to an ``entire'' function $\zeta(x,y)$ on the space
${\mathbb S}_\infty:=\C_\infty^\ast \times \Zp$, where $\C_\infty$ was defined in Section \ref{view} above.
In practical terms this means that for each $y\in \Zp$ we obtain an entire power series $\zeta(x,y)$
in $x^{-1}$ {\it with coefficients in $K$}.

\begin{theorem}\label{swan1}
Fix $y\in \Zp$. Then each zero of $\zeta(x,y)$ is uniquely determined by its
absolute value (including multiplicity).
\end{theorem}

\noindent
Theorem \ref{swan1} is the result of Wan ($q=p$) \cite{wa1}, Poonen, and Sheats (arbitrary $q$) \cite{sh1}, and
the next result then follows immediately.
\begin{cor}\label{swan2} The zeroes of $\zeta(s)$ lie in $K$ and are simple.
\end{cor}

We view Theorem \ref{swan1} as the ``Riemann Hypothesis'' for $\zeta(s)$. Obviously
as elements of $K$, the zeroes of $\zeta(s)$ are trivially unramified over $K$.
\subsection{The failure of the naive generalization of Wan/Sheats}\label{failure}
Let $L/k$ be a finite abelian extension of order $\beta$ prime to $p$ and exponent greater
than $p$; such extensions are very easy to construct. Let  $G:={\rm Gal}(L/k)$
and let $\mathfrak O$ be the ring of $A$-integers of $L$.

Using the monic generators of the ideal norms from $\mathfrak O$ to $A$, one is
easily able to define the analog $\zeta_\mathfrak O(s)$ of the {\it Dedekind zeta function} of $\mathfrak O$.
It is easy to see that $\zeta_\mathfrak O(x,y)$  also analytically continues to a family of entire functions with
coefficients in the local field $K$.

Let $\bar K$ be an algebraic closure of $K$ as in Section \ref{view} and let
$\psi\colon G\to \bar K^\ast$ be a character. As $[L\colon k]$ is prime to $p$, it
does not matter whether we view $\psi$ as having values in $\bar K$ or in $\C$.
Using the formalism of the classical zeta factorizations (i.e., we use the same local
polynomials  but then reduce them modulo $p$), one obtains the expected factorization
\begin{equation}\label{fail1}
\zeta_\mathfrak O(s)=\prod_\psi L(\psi,s)\,,
\end{equation}
where $L(\psi,s)$ has the obvious definition. The functions $L(\psi;x,y)$
 also have analytic
continuations as entire power series in analogy with $\zeta(s)$.\medskip 

\noindent
{\rm N.B.:} Let $\zeta_\beta$ be a primitive $\beta-$th root of unity and let
$K_\beta=K(\zeta_\beta)$ be the unramified extension of $K$ obtained by adjoining
$\zeta_\beta$. One sees that $L(\psi,s)$ gives rise to a family of entire power series
in $x^{-1}$ with coefficients in $K_\beta$. In fact, one can readily find examples where
these coefficients do {\em not} lie in $K$ itself.\medskip

As $[K_\beta \colon K]<\infty$, the field $K_\beta$ is automatically complete. On the other
hand, for fixed $y\in \Zp$, the entire function $L(\psi;x,y)$ factors over its zeroes.
If all such zeroes were contained in $K$, then the coefficients must also be in $K$ and we
have seen this is impossible in general.
\subsection{The ``unramified Riemann hypothesis'' in 
finite characteristic}\label{urh}
Continuing with the above discussion, the best that could be hoped for in terms of the
zeroes of $L(\psi; x,y)$ is that they lie in the field $K_\beta$ and so in particular
are unramified over $K$. We state this as a question.
\begin{question}\label{urh1}
Are the zeroes of $L(\psi,s)$ (and thus also $\zeta_\mathfrak O (s)$) unramified over
$K$?
\end{question}

We view the lack of ramification of the zeroes as the appropriate analog of the Riemann hypothesis for
such finite characteristic valued functions.
\begin{rem}\label{urh2}
A completely similar analysis works when $L/k$ is an arbitrary, not necessarily abelian, 
Galois extension of finite order $\beta$ prime to $p$.
We will discuss the case where the order is divisible by $p$ in Subsection \ref{cav} below.
\end{rem}

Merely stipulating
that the zeroes are unramified does not guarantee that they lie in a finite extension.
However, the next question is very natural in light of Wan/Sheats.
\begin{question}\label{urh3}
Let $L/k$ be finite Galois of order $\beta$ prime to $p$. Let $\rho$ be an irreducible
representation of ${\rm Gal}(L/k)$ and let $L(\rho;x,y)$ be the associated
$L$-series (which we also know to be entire, etc.). Fix $y\in\Zp$. Is every zero of
$L(\rho;x,y)$ uniquely determined by its absolute value (including multiplicity)?
\end{question}

Note of course that if the answer to Question \ref{urh3} is yes, then the zeroes of the $L$-series
will be unramified, simple and lie in $K(\zeta_\beta)$, where $\beta=[L\colon k]$, and so are then obviously contained in
a finite extension of $K$.
\begin{rem}\label{urh4}
The simplicity and lack of ramification of such zeroes is also
expected in the analogous classical case. \end{rem}
\begin{rem}\label{urh5}
The factorizations into $L$-series, combined with positive answers to the above
questions, would allow one to predict the fields obtained by adjoining the zeroes of
the analogs of Dedekind zeta functions, etc.\end{rem}
\section{Complements and Caveats}\label{comp}
\subsection{The connection with the group $S_{(q)}$}\label{con}
Let $q$ continue as above and let $y\in \Zp$. Write $y$ $q$-adically 
as
\begin{equation}\label{con1}
y=\sum_{i=0}^\infty c_i q^i
\end{equation}
where $0\leq c_i<q$ for all $i$. 

Now let $\rho$ be a permutation of the set $\{0,1,2,\ldots\}$.

\begin{defn}\label{con2}
We define $\rho_\ast (y)$, $y\in \Zp$, by
\begin{equation}\label{pi(n)2} 
\rho_\ast(y):=\sum_{i=0}^\infty c_i q^{\rho(i)}\,,
\end{equation} where $y$ is written as in Equation \ref{con1}.\end{defn} 

Clearly $y\mapsto \rho_\ast(y)$ gives a representation of $\rho$ as a bijection
of $\Zp$ which is denoted $\rho_\ast$. The group of such bijections is 
is denoted $S_{(q)}$.
In \cite{go4}, it is shown that $\rho_\ast$ is, in fact, a homeomorphism of $\Zp$ which stabilizes
both the nonnegative and nonpositive integers lying inside $\Zp$. Because of this
stability,
one is able to also define a corresponding action of $\rho_\ast$ on $K$, as
certain {\it $\Fq$-linear} continuous operators. In fact,
$\rho_\ast$  appears to also act on the zeta zeroes (though it must be admitted that there
is so far only a bit of evidence for this).

It is by no means clear how, or even if, one can extend $\rho_\ast$ to the algebraic
closure $\bar K$. However, extending it linearity to $\bar{\F}_q$ gives the obvious natural extension to the unramified
closure $K^{\rm un}$ of $K$, and thus fits well with the ideas of Subsection \ref{urh}.

\subsection{Caveats}\label{cav} 
In this final subsection we present some caveats related to Questions \ref{urh1}
and \ref{urh3}. First of all let $L/k$ now be a finite abelian extension of order
$p$ with group $G$. Let $\mathfrak O\subset L$ be the $A$-integers as before.
We assume, for simplicity, that $L/k$ is ramified at one finite prime
$\mathfrak p=(f)$, $f$ monic -- such examples are easy to find. Let
$\zeta_\mathfrak O(s)$ be the associated analog of the Dedekind zeta function.
Using the characteristic $0$ factors as before we still have the factorization
$$\zeta_\mathfrak O(s)=\prod_\psi L(\psi,s).\,$$
where $\psi$ runs through the characteristic $0$ valued characters of $G$
and where we then reduce the local factors modulo $p$. As the order of $G$ is $p$, we then see
immediately that if $\psi$ is nontrivial, $L(\psi,s)=\zeta_A(s)(1-f^{-s})$ giving the
factorization
$$\zeta_\mathfrak O(s)=\zeta_A(s)^p(1-f^{-s})^{p-1}\,.$$
Obviously 
the zeroes of the Euler factor $(1-f^{-s})$ may even be inseparable depending
on the degree of $f$. In these sort of cases, though, it may be possible to systematically
remove these spurious zeroes through a closer examination of ramification.

On the other hand, we should expect all of this to work for Drinfeld's general
rings $A$ (= the affine algebra of a smooth projective geometrically connected
curve over $\Fq$ minus a chosen place $\infty$). However, in this generality, G. B\"ockle
has examples where there are such spurious zeroes which are {\it not}
attributable to issues involving Euler factors. Different techniques will thus be needed here to
handle these.

Moreover, one would expect a similar theory to hold at the interpolations of these
functions at finite primes $f$. Here one has to consistently remove Euler factors
of the form $(1-f^ix^{-d})$ where $d$ is the degree of $f$ in order to even do
the interpolations. Again removing such Euler factors will contribute many spurious zeroes.

Thus, until a theory arises to handle the spurious zeroes, the general versions
of Questions \ref{urh1} and \ref{urh3} -- applying to all rings $A$ and all places of the
associated complete curve -- should be modified to pertain to
{\em almost all} zeroes for a given $y\in \Zp$, etc.

\end{document}